\documentclass[12pt]{article}                                               
\input{amssym.def}                                                           
\input{amssym.tex}                                                           
                                                                             
\begin{document}                                                             
%**********************************************                              
\title{Mapping class groups are linear}

\author{Igor Nikolaev
\footnote{Partially supported 
by NSERC.}}

%**************************************************
\date{}
 \maketitle
%**************************************************

\newtheorem{thm}{Theorem}
\newtheorem{lem}{Lemma}
\newtheorem{dfn}{Definition}
\newtheorem{rmk}{Remark}
\newtheorem{cor}{Corollary}
\newtheorem{prp}{Proposition}
\newtheorem{exm}{Example}
%*************************************************
%******************************************************************
\begin{abstract}
It is shown,  that the mapping class group of a surface of genus $g\ge 2$ 
admits a faithful representation into the  matrix group $GL_{6g-6}({\Bbb Z})$.  
The proof is based on a categorical correspondence between the Riemann surfaces
and the so-called toric $AF$-algebras.

\vspace{7mm}


\vspace{5mm}
{\it AMS (MOS) Subj. Class.:  14H55, 46L85}
\end{abstract}

%**************************************************************************
\section{Introduction}
%***************************************************************************
{\bf A. The Harvey conjecture.}
The mapping class group has  been introduced in the 1920-ies by M.~Dehn \cite{Deh1}. 
Such a  group, $Mod~(X)$, is defined as the group of isotopy classes of the
orientation-preserving diffeomorphisms of a two-sided closed surface $X$ of
genus $g\ge 1$.  The group is known to be  prominent in  algebraic geometry \cite{HaLo1}, 
topology \cite{Thu1}  and dynamics \cite{Thu2}.  When $X$
is a   torus, the $Mod~(X)$ is isomorphic to the group  $SL_2({\Bbb Z})$. 
(The $SL_2({\Bbb Z})$ is called a modular group, hence our notation for the mapping class group.) 
A little is known about the representations of $Mod~(X)$ beyond the case $g=1$.
Recall, that the group is called {\it linear},  if there exists a faithful
representation into the matrix group $GL_m(R)$, where $R$ is a commutative ring.
The braid groups are known  to be linear \cite{Big1}.  Using a modification of the argument for 
the braid groups, it is possible to prove, that $Mod~(X)$ is linear in the case $g=2$
\cite{BiBu1}. Whether the mapping class group is linear for  $g\ge 3$, 
is an open problem, known as the {\it Harvey conjecture} \cite{Har1}, p.267.

\medskip\noindent
{\bf B. The Teichmueller functor.}
A covariant (non-injective)  functor from a category of generic Riemann
surfaces to a category of the so-called toric $AF$-algebras 
(see section ~2 for a definition) was 
constructed in \cite{Nik1}. The functor, a {\it Teichmueller functor},  
maps any  pair of isomorphic (i.e.  conformal) Riemann surfaces to a pair of 
the stably isomorphic (Morita equivalent) toric  $AF$-algebras.  
Since each isomorphism of Riemann surfaces is given by an element 
of $Mod~(X)$ \cite{HaLo1}, it is natural to ask about a representation of $Mod~(X)$
by the stable isomorphisms of  toric $AF$-algebras; roughly, our objective can be stated 
as follows.

\medskip\noindent
{\sf Main problem.} 
{\it To study the Harvey conjecture from the standpoint of 
  toric $AF$-algebras.}

\medskip\noindent
The stable isomorphisms of toric $AF$-algebras are well understood and
surprisingly simple;  provided  the automorphism group of
the algebra is trivial (this is true for a generic algebra),  its group of stable isomorphism  
admits a faithful representation into a matrix group over the commutative ring ${\Bbb Z}$ \cite{E}.  
This fact, combined  with the properties of the  Teichmueller functor, implies an amazingly simple positive
 solution to the Harvey conjecture.   
%*******************************************************************
\begin{thm}\label{thm1}
For every surface $X$ of genus $g\ge 2$,  there exists a faithful representation 
$\rho: Mod~(X)\rightarrow GL_{6g-6}({\Bbb Z})$.
\end{thm}
%******************************************************************* 
The structure of the note is as follows. In section 2, 
the preliminary facts, necessary to prove theorem \ref{thm1},  
are brought together. Theorem \ref{thm1} is proved in section 3.

%**************************************************************************
\section{Preliminaries}
%***************************************************************************
We review the toric $AF$-algebras and the Teichmueller functor on the space of generic Riemann surfaces;
the reader is referred to \cite{B}, \cite{Bra1},  \cite{E} and \cite{Nik1} for details.

%**************************************************************************
\subsection{$AF$-algebras}
%***************************************************************************
{\bf A. The $C^*$-algebras.}
By the $C^*$-algebra one understands a noncommutative Banach
algebra with an involution. Namely, a $C^*$-algebra
$A$ is an algebra over the complex numbers $\Bbb C$ with a norm $a\mapsto ||a||$
and an involution $a\mapsto a^*, a\in A$, such that $A$ is
complete with the respect to the norm, and such that
$||ab||\le ||a||~||b||$ and $||a^*a||=||a||^2$ for every
$a,b\in A$. If  $A$ is commutative, then the Gelfand
theorem says that $A$ is isometrically $*$-isomorphic
to the $C^*$-algebra $C_0(X)$ of the continuous complex-valued
functions on a locally compact Hausdorff space $X$.
For otherwise, the algebra $A$ represents a noncommutative topological
space.

\medskip\noindent
{\bf B.  The stable isomorphisms of $C^*$-algebras.} 
Let $A$ be a $C^*$-algebra deemed as a noncommutative
topological space. One can ask, when two such topological spaces
$A,A'$ are homeomorphic? To answer the question, let us recall
the topological $K$-theory.  If $X$ is a (commutative) topological
space, denote by $V_{\Bbb C}(X)$ an abelian monoid consisting
of the isomorphism classes of the complex vector bundles over $X$
endowed with the Whitney sum. The abelian monoid $V_{\Bbb C}(X)$
can be made to an abelian group, $K(X)$, using the Grothendieck
completion. The covariant functor $F: X\to K(X)$ is known to
map the homeomorphic topological spaces $X,X'$ to the isomorphic
abelian groups $K(X), K(X')$. Let  $A,A'$ be the $C^*$-algebras. If one wishes to
define a homeomorphism between the noncommutative topological spaces $A$ and $A'$, 
it will suffice to define an isomorphism between the abelian monoids $V_{\Bbb C}(A)$
and $V_{\Bbb C}(A')$ as suggested by the topological $K$-theory. 
The role of the complex vector bundle of the degree $n$ over the 
$C^*$-algebra $A$ is played by a $C^*$-algebra $M_n(A)=A\otimes M_n$,
i.e. the matrix algebra with the entries in $A$.  The abelian monoid
$V_{\Bbb C}(A)=\cup_{n=1}^{\infty} M_n(A)$ replaces the monoid 
$V_{\Bbb C}(X)$ of the topological $K$-theory. Therefore, 
the noncommutative topological spaces $A,A'$ are homeomorphic,
if the abelian monoids $V_{\Bbb C}(A)\cong V_{\Bbb C}(A')$ are isomorphic. 
The latter equivalence is called a {\it stable isomorphism}
of the $C^*$-algebras $A$ and $A'$ and is formally written as 
$A\otimes {\goth K}\cong A'\otimes {\goth K}$, where 
${\goth K}=\cup_{n=1}^{\infty}M_n$ is the $C^*$-algebra of compact
operators.  Roughly speaking, the stable isomorphism between
the $C^*$-algebras  means that they  are
homeomorphic as the noncommutative topological spaces.

\medskip\noindent
{\bf C. The $AF$-algebras.}
An {\it $AF$-algebra}  (approximately finite $C^*$-algebra) is defined to
be the  norm closure of an ascending sequence of the finite dimensional
$C^*$-algebras $M_n$'s, where  $M_n$ is the $C^*$-algebra of the $n\times n$ matrices
with the entries in ${\Bbb C}$. Here the index $n=(n_1,\dots,n_k)$ represents
a semi-simple matrix algebra $M_n=M_{n_1}\oplus\dots\oplus M_{n_k}$.
The ascending sequence mentioned above  can be written as 
%***********************************************************
$M_1\buildrel\rm\varphi_1\over\longrightarrow M_2
   \buildrel\rm\varphi_2\over\longrightarrow\dots,
$
%****************************************************
where $M_i$ are the finite dimensional $C^*$-algebras and
$\varphi_i$ the homomorphisms between such algebras.  The set-theoretic limit
$A=\lim M_n$ has a natural algebraic structure given by the formula
$a_m+b_k\to a+b$; here $a_m\to a,b_k\to b$ for the
sequences $a_m\in M_m,b_k\in M_k$.  
The homomorphisms $\varphi_i$ can be arranged into  a graph as follows. 
Let  $M_i=M_{i_1}\oplus\dots\oplus M_{i_k}$ and 
$M_{i'}=M_{i_1'}\oplus\dots\oplus M_{i_k'}$ be 
the semi-simple $C^*$-algebras and $\varphi_i: M_i\to M_{i'}$ the  homomorphism. 
One has the two sets of vertices $V_{i_1},\dots, V_{i_k}$ and $V_{i_1'},\dots, V_{i_k'}$
joined by the $a_{rs}$ edges, whenever the summand $M_{i_r}$ contains $a_{rs}$
copies of the summand $M_{i_s'}$ under the embedding $\varphi_i$. 
As $i$ varies, one obtains an infinite graph called a {\it Bratteli diagram} of the
$AF$-algebra \cite{Bra1}. The Bratteli diagram defines a unique  $AF$-algebra.

\medskip\noindent
{\bf D. The stationary $AF$-algebras.}
If the homomorphisms $\varphi_1 =\varphi_2=\dots=Const$ in the definition of 
the $AF$-algebra $A$,  the $AF$-algebra $A$ is called {\it stationary}. 
The Bratteli diagram of a stationary $AF$-algebra looks like a periodic 
graph  with the incidence matrix $A=(a_{rs})$ repeated over and over again. 
Since  matrix $A$ is a non-negative integer matrix, one can take a power of
$A$ to obtain a strictly positive integer matrix -- which we always assume 
to be the case.  The stationary $AF$-algebra has a non-trivial group
of the automorphisms \cite{E}, Ch.6.

%**************************************************************************
\subsection{Teichmueller functor}
%***************************************************************************
{\bf A. The Jacobi-Perron continued fractions.}
Let $\lambda=(\lambda_1,\dots,\lambda_n)$ be a vector with the non-negative
real entries, such that $\lambda_1\ne 0$;  consider a projective class 
$(1,\theta_1,\dots,\theta_{n-1})$ of $\lambda$, where $\theta_{i-1}={\lambda_i\over\lambda_1}$
for $1\le i\le n$.  The continued fraction
%***********************************************************************
%\footnote{The reader can verify, that the object is indeed a continued fraction,
%which generalizes the regular continued fraction $\theta= [b_1, b_2,\dots]$,
%written in the matrix form 
%$\left(\small\matrix{1\cr \theta}\right)=
%\lim_{k\to\infty} \left(\small\matrix{0 & 1\cr 1 & b_1}\right)\dots
%\left(\small\matrix{0 & 1\cr 1 & b_k}\right)
%\left(\small\matrix{0\cr 1}\right).$}
%**************************************************************
%\begin{equation}
$$
\left(\matrix{1\cr \theta_1\cr\vdots\cr\theta_{n-1}} \right)=
\lim_{k\to\infty} 
\left(\matrix{0 &  0 & \dots & 0 & 1\cr
              1 &  0 & \dots & 0 & b_1^{(1)}\cr
              \vdots &\vdots & &\vdots &\vdots\cr
              0 &  0 & \dots & 1 & b_{n-1}^{(1)}}\right)
\dots 
\left(\matrix{0 &  0 & \dots & 0 & 1\cr
              1 &  0 & \dots & 0 & b_1^{(k)}\cr
              \vdots &\vdots & &\vdots &\vdots\cr
              0 &  0 & \dots & 1 & b_{n-1}^{(k)}}\right)
\left(\matrix{0\cr 0\cr\vdots\cr 1} \right),
$$
%\end{equation}
%**************************************************************
where $b_i^{(j)}\in {\Bbb N}\cup\{0\}$, is called the {\it Jacobi-Perron
fraction}. To recover the integers $b_i^{(k)}$ from the vector $(\theta_1,\dots,\theta_{n-1})$,
one has to repeatedly solve the following system of equations:
$ \theta_1 = b_1^{(1)} +  {1\over\theta_{n-1}'}$,
$\theta_2 = b_2^{(1)} +  {\theta_1'\over\theta_{n-1}'}$,
$\theta_{n-1} = b_{n-1}^{(1)} +  {\theta_{n-2}'\over\theta_{n-1}'}$,
where $(\theta_1',\dots,\theta_{n-1}')$ is the next input vector.
Thus, each vector $(\theta_1,\dots,\theta_{n-1})$ gives rise to a formal 
Jacobi-Perron continued fraction, which can be convergent or divergent.
We let $A^{(0)}=\delta_{ij}$ (the Kronecker delta) and
$
A_i^{(k+n)}=\sum_{j=0}^{n-1} b_i^{(k)}A_i^{(\nu +j)}, \quad b_0^{(k)}=1,
$
where $i=0,\dots,n-1$ and $k=0,1,\dots,\infty$. 
The Jacobi-Perron continued fraction  of vector 
$(\theta_1,\dots,\theta_{n-1})$  is said to be {\it convergent}, if 
$\theta_i=\lim_{k\to\infty}{A_i^{(k)}\over A_0^{(k)}}$ 
for all $i=1,\dots,n-1$.
Unless $n=2$, the convergence of the individual Jacobi-Perron fraction is a 
difficult problem;  however, it is known that the Jacobi-Perron fractions  
converge for a  generic subset of the vectors in the space ${\Bbb R}^{n-1}$
\cite{Bau1}.

\medskip\noindent
{\bf B. The toric $AF$-algebras.}
Denote by $T_S(g)$ the Teichmueller
space of genus $g\ge 1$ with a distinguished point $S$. Let
$q\in H^0(S, \Omega^{\otimes 2})$ be a holomorphic quadratic 
differential on the Riemann surface $S$, such that all zeroes
of $q$ (if any) are simple. By $\widetilde S$ we mean a double 
cover of $S$ ramified over the zeroes of $q$ and by
$H_1^{odd}(\widetilde S)$ the odd part of the integral homology of $\widetilde S$ relative to the  zeroes.
Note that $H_1^{odd}(\widetilde S)\cong {\Bbb Z}^{n}$, where $n=6g-6$ if $g\ge 2$ and $n=2$ if $g=1$. 
It follows from the Main Theorem of \cite{HuMa1}, that
$T_S(g)-\{pt\}\cong Hom~(H_1^{odd}(\widetilde S); {\Bbb R})-\{0\}$, where $0$ is the zero homomorphism
%**********************************************************************************
\footnote{To be precise, the theorem mentions a local homeomorphism
$h: Hom~(H_1^{odd}(\widetilde S); {\Bbb R})-\{0\}\to T_S(g)-\{pt\}$ 
\cite{HuMa1}, p.222. Since $T_S(g)$ is simply connected, $h$ extends to
a  global homeomorphism between the two spaces. Indeed, let
$\lambda\in  Hom~(H_1^{odd}(\widetilde S); {\Bbb R})$. 
It is easy to see, that $\lambda=0$ corresponds to the distinguished
point $S\in T_S(g)$, while $\lambda=\infty$ represent the boundary
of the space $T_S(g)$; thus, every ball $|\lambda|<C$ is homotopy equivalent
 to the ball  $|\lambda|<\infty$. Note also, that we are interested in
$q$'s with the generic (simple) zeroes; the higher order zeroes --
which can be an obstacle in the construction of global coordinates --
are excluded. The interested reader can consult \cite{Thu2}, p.425 (the last
paragraph) for the details; the mentioned there piecewise linear integral structure
breaks $T_S(g)$ into a {\it finite} number of cones issued from $S$ and in 
these terms our construction means that we take a cone and extend it (by the linearity)
to the entire $T_S(g)$. The cones differ from each other by a permutation on
the set $\lambda=(\lambda_1,\dots,\lambda_n)$; distinct permutations correspond
to the different coordinates in the space $T_S(g)$.}
%**********************************************************************************
.  Finally, denote by $\lambda=(\lambda_1,\dots,\lambda_{n})$ the image of a basis of 
$H_1^{odd}(\widetilde S)$ in the real line ${\Bbb R}$, such that $\lambda_1\ne0$. 
(Note that such an option always exists, since the zero homomorphism is excluded.) 
We let $\theta=(\theta_1,\dots,\theta_{n-1})$, where $\theta_i=\lambda_{i-1}/\lambda_1$. 
Recall that,  up to a scalar multiple, the vector $(1,\theta)\in {\Bbb R}^{n}$ is the limit
of a generically convergent Jacobi-Perron continued fraction:
%**************************************************************
$$
\left(\matrix{1\cr \theta}\right)=
\lim_{k\to\infty} \left(\matrix{0 & 1\cr I & b_1}\right)\dots
\left(\matrix{0 & 1\cr I & b_k}\right)
\left(\matrix{0\cr {\Bbb I}}\right),
$$
%****************************************************************
where $b_i=(b^{(i)}_1,\dots, b^{(i)}_{n-1})^T$ is a vector of the non-negative integers,  
$I$ the unit matrix and ${\Bbb I}=(0,\dots, 0, 1)^T$. We introduce an  $AF$-algebra, 
${\Bbb A}_{\theta}$, via the Bratteli diagram, shown in Fig.1.   
(The numbers $b_j^{(i)}$ of the diagram  indicate the multiplicity  of edges of 
the graph.)  Let us call ${\Bbb A}_{\theta}$ a {\it toric $AF$-algebra}. 
Note that in the $g=1$ case, the Jacobi-Perron 
fraction coincides with the  regular continued fraction and ${\Bbb A}_{\theta}$
becomes the Effros-Shen $AF$-algebra of a noncommutative torus
\cite{EfSh1}.

%\bigskip
%*******************************************************************
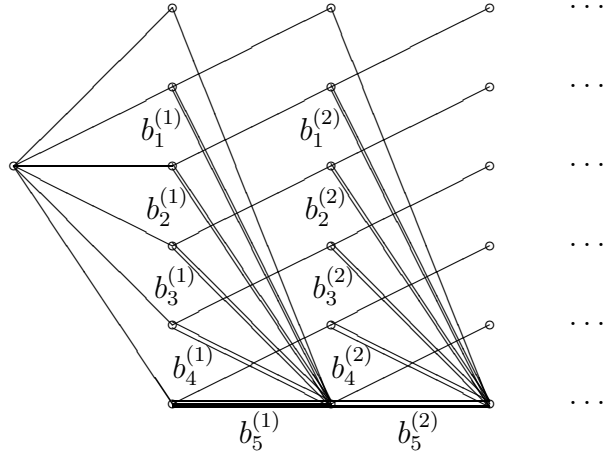
\begin{figure}[here]
\begin{picture}(300,200)(0,0)

\put(40,110){\circle{3}}
%%%%%%%%%%%%%%%%%%%%%%%%%%%%%%%%%%%%

\put(40,110){\line(1,1){60}}
\put(40,110){\line(2,1){60}}
\put(40,110){\line(1,0){60}}
\put(40,110){\line(2,-1){60}}
\put(40,110){\line(1,-1){60}}
\put(40,110){\line(2,-3){60}}
%%%%%%%%%%%%%%%%%%%%%%%%%%%%%%%%%%%%%%%

\put(100,20){\circle{3}}
\put(100,50){\circle{3}}
\put(100,80){\circle{3}}
\put(100,110){\circle{3}}
\put(100,140){\circle{3}}
\put(100,170){\circle{3}}

\put(160,20){\circle{3}}
\put(160,50){\circle{3}}
\put(160,80){\circle{3}}
\put(160,110){\circle{3}}
\put(160,140){\circle{3}}
\put(160,170){\circle{3}}

\put(220,20){\circle{3}}
\put(220,50){\circle{3}}
\put(220,80){\circle{3}}
\put(220,110){\circle{3}}
\put(220,140){\circle{3}}
\put(220,170){\circle{3}}

%%%%%%%%%%%%%%%%%%%%%%%%%%%%%%%%%%%%%%%%%%%%%%%%%%%%%%%%%%%%%%%
\put(160,20){\line(2,1){60}}
\put(160,19){\line(1,0){60}}
%\put(160,20){\line(1,0){60}}
\put(160,21){\line(1,0){60}}
\put(160,50){\line(2,1){60}}
\put(160,49){\line(2,-1){60}}
\put(160,51){\line(2,-1){60}}
\put(160,80){\line(2,1){60}}
\put(160,79){\line(1,-1){60}}
\put(160,81){\line(1,-1){60}}
\put(160,110){\line(2,1){60}}
\put(160,109){\line(2,-3){60}}
\put(160,111){\line(2,-3){60}}
\put(160,140){\line(2,1){60}}
\put(160,139){\line(1,-2){60}}
\put(160,141){\line(1,-2){60}}
\put(160,170){\line(2,-5){60}}

%%%%%%%%%%%%%%%%%%%%%%%%%%%%%%%%%%%

\put(100,20){\line(2,1){60}}
\put(100,19){\line(1,0){60}}
\put(100,20){\line(1,0){60}}
\put(100,21){\line(1,0){60}}
\put(100,50){\line(2,1){60}}
\put(100,49){\line(2,-1){60}}
\put(100,51){\line(2,-1){60}}
\put(100,80){\line(2,1){60}}
\put(100,79){\line(1,-1){60}}
\put(100,81){\line(1,-1){60}}
\put(100,110){\line(2,1){60}}
\put(100,109){\line(2,-3){60}}
\put(100,111){\line(2,-3){60}}
\put(100,140){\line(2,1){60}}
\put(100,139){\line(1,-2){60}}
\put(100,141){\line(1,-2){60}}
\put(100,170){\line(2,-5){60}}

%%%%%%%%%%%%%%%%%%%%%%%%%%%%%%%%%%%

\put(250,20){$\dots$}
\put(250,50){$\dots$}
\put(250,80){$\dots$}
\put(250,110){$\dots$}
\put(250,140){$\dots$}
\put(250,170){$\dots$}

%%%%%%%%%%%%%%%%%%%%%%%%%%%%%%%%%%%%%%%%%%%

\put(125,5){$b_5^{(1)}$}
\put(100,30){$b_4^{(1)}$}
\put(93,60){$b_3^{(1)}$}
\put(90,90){$b_2^{(1)}$}
\put(88,120){$b_1^{(1)}$}

\put(185,5){$b_5^{(2)}$}
\put(160,30){$b_4^{(2)}$}
\put(153,60){$b_3^{(2)}$}
\put(150,90){$b_2^{(2)}$}
\put(148,120){$b_1^{(2)}$}

%%%%%%%%%%%%%%%%%%%%%%%%%%%%%%%%%%%%%%%%%%

\end{picture}

\caption{The Bratteli diagram of a toric $AF$-algebra of genus $2$.}
\end{figure}
%*******************************************************************

\medskip\noindent
{\bf C. The Teichmueller functor.}
Denote by $V$ the maximal subset of $T_S(g)$,  such that for each Riemann surface $R\in V$, 
there exists a convergent Jacobi-Perron continued fraction. 
Let $F$ be the map which sends the Riemann surfaces into the toric $AF$-algebras
according to the formula  $R\mapsto {\Bbb A}_{\theta}$.  
Finally, let $W$ be the image of  $V$ under the mapping $F$.
%****************************************************************
\begin{lem}\label{lm1}
{\bf (\cite{Nik1})}
The set $V$ is a generic subset of $T_S(g)$ and the  map $F$ has the following properties:

\medskip
(i) $V\cong W\times (0,\infty)$ is a trivial fiber bundle, whose
projection  map $p: V\to W$ coincides with $F$;

\smallskip
(ii) $F: V\to W$ is a covariant functor, which maps isomorphic Riemann surfaces 
$R,R'\in V$ to  stably isomorphic toric $AF$-algebras ${\Bbb A}_{\theta},{\Bbb A}_{\theta'}\in W$.
\end{lem}
%***************************************************************

%**************************************************************************
\section{Proof of theorem 1}
%***************************************************************************
As before, let $W$ denote  the set of toric $AF$-algebras of genus $g\ge 2$. 
Let $G$ be a finitely presented group and  $G\times W\to W$ its action on $W$ 
by  the stable isomorphisms of toric $AF$-algebras; in other words, 
 $\gamma ({\Bbb A}_{\theta})\otimes {\goth K}\cong
{\Bbb A}_{\theta}\otimes {\goth K}$  for all $\gamma\in G$ and all 
${\Bbb A}_{\theta}\in W$.   The following preparatory lemma will be 
important.
%**********************************************************************
\begin{lem}\label{lm2}
For each ${\Bbb A}_{\theta}\in W$,  there exists a 
representation $\rho_{{\Bbb A}_{\theta}}: G\to GL_{6g-6}({\Bbb Z})$.  
\end{lem}
%**********************************************************************
{\it Proof.}
The proof of  lemma is based on the following well known criterion
of the stable isomorphism for the (toric)  $AF$-algebras: a pair of such algebras
${\Bbb A}_{\theta}, {\Bbb A}_{\theta'}$ are stably isomorphic if and only
if their Bratteli diagrams coincide, except (possibly) a finite part
of the diagram, see \cite{E}, Theorem 2.3. (Note, that the order isomorphism
between  the dimension groups, mentioned in the original text, translates 
to the language of the Bratteli diagrams as stated.)

Let $G$ be a finitely presented group on the generators  $\{\gamma_1, \dots, \gamma_m\}$
subject to relations $r_1,\dots,r_n$.  Let  ${\Bbb A}_{\theta}\in W$.
Since $G$ acts on the toric $AF$-algebra ${\Bbb A}_{\theta}$ by stable isomorphisms,
the toric $AF$-algebras ${\Bbb A}_{\theta_1}:=\gamma_1({\Bbb A}_{\theta}),\dots,   
{\Bbb A}_{\theta_m}:=\gamma_m({\Bbb A}_{\theta})$ are stably isomorphic to 
${\Bbb A}_{\theta}$; moreover, by transitivity, they are also pairwise stably isomorphic.
Therefore, the Bratteli diagrams of ${\Bbb A}_{\theta_1},\dots, {\Bbb A}_{\theta_m}$ 
 coincide everywhere except, possibly, some finite parts.  
We shall denote by ${\Bbb A}_{\theta_{\max}}\in W$
a toric $AF$-algebra, whose Bratteli diagram is the maximal common part
of the Bratteli diagrams of ${\Bbb A}_{\theta_i}$ for $1\le i\le m$;
such a choice is unique and defined correctly because the set $\{{\Bbb A}_{\theta_i}\}$
is a finite set.  By the definition of a toric $AF$-algebra, the vectors 
$\theta_i=(1,\theta_1^{(i)},\dots,\theta_{6g-7}^{(i)})$
are related to the vector $\theta_{\max}=(1, \theta_1^{(\max)},\dots,\theta_{6g-7}^{(\max)})$ 
by the formula:
%*******************************************************************
$$
\left(\matrix{1\cr \theta_1^{(i)}\cr\vdots\cr\theta_{6g-7}^{(i)}} \right)
=\underbrace{
\left(\matrix{0 &  0 & \dots & 0 & 1\cr
              1 &  0 & \dots & 0 & b_1^{(1)(i)}\cr
              \vdots &\vdots & &\vdots &\vdots\cr
              0 &  0 & \dots & 1 & b_{6g-7}^{(1)(i)}}\right)
\dots 
\left(\matrix{0 &  0 & \dots & 0 & 1\cr
              1 &  0 & \dots & 0 & b_1^{(k)(i)}\cr
              \vdots &\vdots & &\vdots &\vdots\cr
              0 &  0 & \dots & 1 & b_{6g-7}^{(k)(i)}}\right)
}_{A_i}
\left(\matrix{1\cr \theta^{(\max)}_1\cr\vdots\cr\theta^{(\max)}_{6g-7}} \right)
$$
%***************************************************************** 
The above expression can be written in the matrix form $\theta_i=A_i\theta_{\max}$, where 
$A_i\in GL_{6g-6}({\Bbb Z})$. Thus, one gets a matrix representation of the
generator $\gamma_i$,  given by the formula $\rho_{{\Bbb A}_{\theta}}(\gamma_i):=A_i$.

The map $\rho_{{\Bbb A}_{\theta}}: G\to GL_{6g-6}({\Bbb Z})$ extends to the rest of the group $G$
 via its values on the generators;  namely,  for every $g\in G$ one sets $\rho_{{\Bbb A}_{\theta}}(g)= A_1^{k_1}\dots A_m^{k_m}$,
whenever $g=\gamma_1^{k_1}\dots \gamma_m^{k_m}$. Let us verify, that the 
map $\rho_{{\Bbb A}_{\theta}}$ is a well defined homomorphism of groups $G$ and $GL_{6g-6}({\Bbb Z})$.
Indeed, let us write $g_1=\gamma_1^{k_1}\dots\gamma_m^{k_m}$ and
$g_2=\gamma_1^{s_1}\dots\gamma_m^{s_m}$ for a pair of elements $g_1,g_2\in G$;
then their product 
$g_1g_2=\gamma_1^{k_1}\dots\gamma_m^{k_m}\gamma_1^{s_1}\dots\gamma_m^{s_m}=
\gamma_1^{l_1}\dots\gamma_m^{l_m}$,
where the last equality is obtained by a reduction of words using the 
relations $r_1,\dots,r_n$.  One can write relations $r_i$ in their matrix
form $\rho_{{\Bbb A}_{\theta}}(r_i)$; thus, one gets the matrix equality
$A_1^{l_1}\dots A_m^{l_m}= A_1^{k_1}\dots A_m^{k_m}A_1^{s_1}\dots A_m^{s_m}$.
It is immediate from the last  equation, that 
$\rho_{{\Bbb A}_{\theta}}(g_1g_2)=
A_1^{l_1}\dots A_m^{l_m}= A_1^{k_1}\dots A_m^{k_m}A_1^{s_1}\dots A_m^{s_m}=
\rho_{{\Bbb A}_{\theta}}(g_1)\rho_{{\Bbb A}_{\theta}}(g_2)$ for  $\forall g_1,g_2\in G$,
i.e.  $\rho_{{\Bbb A}_{\theta}}$ is a homomorphism. Lemma \ref{lm2} follows.
$\square$

\bigskip
Let  $W_{aper}\subset W$ be a set consisting  of the toric $AF$-algebras, whose  Bratteli diagrams
do not contain periodic (infinitely repeated) blocks;  these are known as  non-stationary toric  $AF$-algebras
and they are generic in the set $W$ endowed with the natural topology. 
We call  the action of $G$  on the toric $AF$-algebra  ${\Bbb A}_{\theta}\in W$
{\it free}, if   $\gamma ({\Bbb A}_{\theta})={\Bbb A}_{\theta}$ implies $\gamma=Id$.
%**********************************************************************
\begin{lem}\label{lm3}
Let ${\Bbb A}_{\theta}\in W_{aper}$  and  $G$ be  
free on the ${\Bbb A}_{\theta}$.  Then  $\rho_{{\Bbb A}_{\theta}}$ is a 
faithful representation.
\end{lem}
%**********************************************************************
{\it Proof.}
Since the action of $G$ is free, to prove that  $\rho_{{\Bbb A}_{\theta}}$  is faithful, 
it remains  to show, that in the formula  $\theta_i=A_i\theta_{\max}$, it holds
$A_i=I$, if and only if,  $\theta_i=\theta_{\max}$, where $I$
is the unit matrix.   Indeed,  it is immediate that $A_i=I$ implies $\theta_i=\theta_{\max}$.
Suppose now that  $\theta_i=\theta_{\max}$ and, let to the contrary, $A_i\ne I$. 
One gets $\theta_i=A_i \theta_{\max}=\theta_{\max}$.  Such an equation has  a non-trivial solution, 
if and only if,  the vector $\theta_{\max}$ has a periodic 
Jacobi-Perron fraction;  the period of  such a fraction  is given by the matrix $A_i$. This 
is impossible, since it has been assumed, that ${\Bbb A}_{\theta_{\max}}\in W_{aper}$.
The contradiction finishes the proof of lemma \ref{lm3}.
$\square$

\medskip
Let $G=Mod~(X)$, where $X$ is a surface of genus $g\ge 2$. The group $G$
is finitely presented \cite{Deh1}; it  acts on the Teichmueller space $T(g)$
by isomorphisms of the Riemann surfaces. Moreover, the action of $G$ is free on a 
generic set, $U\subset T(g)$,  consisting of the Riemann surfaces with the trivial group 
of the automorphisms. Recall,  that there exists a functor  $F: V\to W$  between the Riemann surfaces and 
toric $AF$-algebras, see lemma  \ref{lm1}.
%**********************************************************************
\begin{lem}\label{lm4}
The pre-image  $F^{-1}(W_{aper})$ is a generic set in the space  $T(g)$.
\end{lem}
%**********************************************************************
{\it Proof.}
Note,  that the set of stationary toric $AF$-algebras is a countable set.
The functor $F$ is a surjective map,  which is continuous with respect to the natural topology 
on the sets $V$ and $W$.  Therefore,  the pre-image of the complement of a countable set is a generic set. 
$\square$

\medskip
To finish the proof,  consider the set $U\cap F^{-1}(W_{aper})$;
this is a non-empty set, since it is the intersection of  two generic subsets of $T(g)$.
Let $R$ be a point (a Riemann surface) in the above set.  In  view of lemma \ref{lm1},
group $G$ acts on the toric $AF$-algebra ${\Bbb A}_{\theta}=F(R)$ by the stable
isomorphisms.  By the construction, the action is free and ${\Bbb A}_{\theta}\in W_{aper}$.
In view of lemma \ref{lm3},  one gets a faithful representation $\rho=\rho_{{\Bbb A}_{\theta}}$
of the group $G=Mod~(X)$ into the matrix group $GL_{6g-6}({\Bbb Z})$. Theorem \ref{thm1}
follows.
$\square$

%**************************************************************************

%**********************************************************

\vskip1cm

\textsc{The Fields Institute for Mathematical Sciences, Toronto, ON, Canada,  
E-mail:} {\sf igor.v.nikolaev@gmail.com}

\smallskip
{\it Current address: 101-315 Holmwood Ave., Ottawa, ON, Canada, K1S 2R2}  

\end{document}